\ifx\shlhetal\undefinedcontrolsequence\let\shlhetal\relax\fi

\documentclass[11pt]{amsart}

\usepackage{amsmath}
\usepackage{amssymb}

\newtheorem{theorem}{Theorem}[section]
\newtheorem{claim}[theorem]{Claim}

\newtheorem{proposition}[theorem]{Proposition}
\newtheorem{corollary}[theorem]{Corollary}

\theoremstyle{definition}
\newtheorem{definition}[theorem]{Definition}

\newtheorem{problem}{Problem}[section]
\newtheorem{question}[theorem]{Question}

\theoremstyle{remark}
\newtheorem{remark}[theorem]{Remark}

\newtheorem{conclusion}[theorem]{Conclusion}

\newcount\skewfactor
\def\mathunderaccent#1#2 {\let\theaccent#1\skewfactor#2
\mathpalette\putaccentunder}
\def\putaccentunder#1#2{\oalign{$#1#2$\crcr\hidewidth
\vbox to.2ex{\hbox{$#1\skew\skewfactor\theaccent{}$}\vss}\hidewidth}}

\def\smallbox#1{\leavevmode\thinspace\hbox{\vrule\vtop{\vbox
   {\hrule\kern1pt\hbox{\vphantom{\tt/}\thinspace{\tt#1}\thinspace}}
   \kern1pt\hrule}\vrule}\thinspace}

\newcommand{\cf}{{\rm cf}}

\newcommand{\then}{{\underline{then}}}

\def\qedref#1{$\qed_{\reforiginal{#1}}$}

\title{Partition calculus and cardinal invariants}
\author{Shimon Garti}
\address{Institute of Mathematics
 The Hebrew University of Jerusalem
 Jerusalem 91904, Israel}
\email{shimon.garty@mail.huji.ac.il}

\author{Saharon Shelah}
\address{Institute of Mathematics
 The Hebrew University of Jerusalem
 Jerusalem 91904, Israel
 and  Department of Mathematics
 Rutgers University
 New Brunswick, NJ 08854, USA}
\email{shelah@math.huji.ac.il}
\urladdr{http://www.math.rutgers.edu/\char`\~shelah}

\thanks{First typed: December 2011 \newline Research supported by the United States-Israel Binational Science Foundation. Publication 995 of the second author}
\subjclass[2010] {03E02, 03E04, 03E05, 03E17, 03E35}
\keywords{Partition calculus, Cardinal characteristics}

\begin{document}
\let\labeloriginal\label
\let\reforiginal\ref

\begin{abstract}
We prove that the strong polarized relation $\binom{\theta}{\omega} \rightarrow \binom{\theta}{\omega}^{1,1}_2$, applied simultaneously for every $\theta\in[\aleph_1,2^{\aleph_0}]$, is consistent with ZFC. Consequently, $\binom{inv}{\omega} \rightarrow \binom{inv}{\omega}^{1,1}_2$ is consistent for every cardinal invariant of the continuum. Some results in this direction are generalized to higher cardinals. \newline 
Nous prouvons que la relation polaris\'ee forte $\binom{\theta}{\omega} \rightarrow \binom{\theta}{\omega}^{1,1}_2$, appliqu\'ee simultan\'ement \`a chaque cardinal $\theta\in[\aleph_1,2^{\aleph_0}]$, est en accord avec ZFC. Par cons\'equent, la relation $\binom{inv}{\omega} \rightarrow \binom{inv}{\omega}^{1,1}_2$ est en accord avec ZFC pour chaque caract\'eristique sur le continu. Nous \'etudions plusieurs 
g\'en\'eralisations pour certains cardinaux \'elev\'es.
\end{abstract}

\maketitle

\newpage

\section{introduction}

The strong polarized relation $\binom{\lambda}{\kappa} \rightarrow \binom{\lambda}{\kappa}^{1,1}_2$ means that for every function $c : \lambda \times \kappa \rightarrow 2$ there are $A \subseteq \lambda$ and $B \subseteq \kappa$ such that $|A|=\lambda, |B|=\kappa$ and $c \upharpoonright (A \times B)$ is constant. The history of this relation begins with \cite{MR0081864}, and later \cite{MR0202613}. A comprehensive discussion on the basic results for this relation appears in \cite{williams}. For a modern discussion see \cite{partitions}.

Cardinal invariants of the continuum are discussed in \cite{MR2768685}. Every cardinal invariant isolates some property of the continuum (i.e., $^\omega 2,^\omega \omega$, or $[\omega]^\omega$ and so forth) and seeks for the minimal cardinality of a set with this property. The value of each cardinal invariant belongs to the interval $[\aleph_1,\mathfrak{c}]$, and except of the trivial invariants (which are the first uncountable cardinal, and $\mathfrak{c}$), the value of each invariant can fall on a large spectrum of cardinals in this interval.
We are interested in the following general problem, from \cite{GaThesis}:

\begin{problem}
\label{iinv} 
Cardinal invariants and the polarized relation. \newline
Let \emph{inv} be a cardinal invariant of the continuum. \newline
Is the relation $\binom{inv}{\omega} \rightarrow \binom{inv}{\omega}^{1,1}_2$ consistent with ZFC?
\end{problem}

Since the continuum hypothesis implies $\binom{\aleph_1}{\aleph_0} \nrightarrow \binom{\aleph_1}{\aleph_0}^{1,1}_2$ (as proved in \cite{MR0081864}), and $inv=\aleph_1$ for every cardinal invariant under the continuum hypothesis, we know that the negative relation $\binom{inv}{\omega} \nrightarrow \binom{inv}{\omega}^{1,1}_2$ is always consistent. This is the background behind problem \ref{iinv}.

In \cite{MR2927607} it is proved that $\binom{\mathfrak{c}}{\omega} \rightarrow \binom{\mathfrak{c}}{\omega}^{1,1}_2$ is consistent with ZFC, and one can judge $\mathfrak{c}$ as a cardinal invariant, giving a positive answer (in this case) for the above problem. But in the model constructed in \cite{MR2927607} there exists an uncountable cardinal $\theta<\mathfrak{c}$ so that $\binom{\theta}{\omega} \nrightarrow \binom{\theta}{\omega}^{1,1}_2$. This gives rise to the following:

\begin{problem}
\label{aaalll}
Simultaneous positive relations. \newline 
Is the relation $\binom{\theta}{\omega} \rightarrow \binom{\theta}{\omega}^{1,1}_2$ consistent with ZFC for every $\theta\in[\aleph_1,2^{\aleph_0}]$ simultaneously?
\end{problem}

By the way we mention that the opposite situation holds in the Cohen model. Namely, adding $\lambda$-many Cohen reals implies $\binom{\theta}{\omega} \nrightarrow \binom{\theta}{\omega}^{1,1}_2$ for every $\theta\in[\aleph_1,2^{\aleph_0}]$. An explicit proof can be found in \cite{gash962}, remark 1.4.

Let us state the known results so far. By \cite{MR2927607}, if $\kappa < \mathfrak{s}$ then $\binom{\kappa}{\omega} \rightarrow \binom{\kappa}{\omega}^{1,1}_2$ iff $\cf(\kappa)>\aleph_0$. Hence forcing $\aleph_0<\cf(inv)\leq inv<\mathfrak{s}$ settles problem \ref{iinv} for such an invariant (the cofinality requirement is easy, in general). For instance, it gives the consistency of $\binom{\mathfrak{b}}{\omega} \rightarrow \binom{\mathfrak{b}}{\omega}^{1,1}_2$, as well as $\binom{\mathfrak{a}}{\omega} \rightarrow \binom{\mathfrak{a}}{\omega}^{1,1}_2$, due to \cite{MR1623206} (chapter VI, \S 6).

So we focus on invariants above $\mathfrak{s}$. In a sense, $\mathfrak{s}$ is a natural invariant for getting `downward positive relations' like $\binom{\kappa}{\omega} \rightarrow \binom{\kappa}{\omega}^{1,1}_2$, whenever $\kappa<\mathfrak{s}$. Here we shall see that the reaping number $\mathfrak{r}$ is a natural invariant for `upward positive relations', namely $\binom{\kappa}{\omega} \rightarrow \binom{\kappa}{\omega}^{1,1}_2$ for every $\kappa>\mathfrak{r}$ whose cofinality is large enough. 

Inasmuch as $\mathfrak{r}<\mathfrak{s}$ is consistent with ZFC, we can cover simultaneously every $\theta\in[\aleph_1,2^{\aleph_0}]$. In the model of \cite{MR879489}, $\aleph_1=\mathfrak{r}<\mathfrak{s}=\aleph_2=\mathfrak{c}$. This gives a positive answer to problem \ref{aaalll}, hence also to problem \ref{iinv}, since every cardinal invariant falls into $\{\aleph_1,\aleph_2\}$ in this model. 

Another result is related to $\mathfrak{d}$. It is not known, yet, if one can increase the continuum above $\aleph_2$ while keeping $\mathfrak{r}<\mathfrak{s}$. Anyhow, dealing with the dominating number $\mathfrak{d}$ one can force $\mathfrak{r}<\mathfrak{d}$ for every prescribed regular value of $\mathfrak{d}$ above $\aleph_1$, as proved in \cite{MR1005010}. Consequently, the relation $\binom{\mathfrak{d}}{\omega} \rightarrow \binom{\mathfrak{d}}{\omega}^{1,1}_2$ is consistent with ZFC for arbitrarily large $\mathfrak{d}$. 

Can we generalize these results to uncountable cardinals? We need some large cardinal assumptions. If $\lambda$ is a supercompact cardinal we can force $\mathfrak{r}_\lambda=\mathfrak{u}_\lambda=\lambda^+$, yielding positive relation for every regular cardinal above $\lambda^+$. We believe that some sort of large cardinals assumption is needed, yet supercompactness is not vital. We still do not know what happens in the general case of an uncountable $\lambda$.
If $\mu$ is a singular cardinal (a limit of strongly inaccessibles, or a parallel assumption) then we can increase $2^\mu$ and prove $\binom{\theta}{\mu} \rightarrow\binom{\theta}{\mu}^{1,1}_2$ for many $\theta$-s in the interval $(\mu,2^\mu]$.

We use standard notation. We employ the letters $\theta, \kappa, \lambda, \mu, \chi$ for infinite cardinals, and $\alpha, \beta, \gamma, \delta, \varepsilon, \zeta$ for ordinals. Topological cardinal invariants of the continuum are denoted as in \cite{MR776622} and \cite{MR2768685}. We denote the continuum by $\mathfrak{c}$. For $A,B\subseteq\lambda$ we denote almost inclusion by $\subseteq^*$, so $A\subseteq^*B$ means $|A\setminus B|<\lambda$. For a regular cardinal $\kappa$ we denote the ideal of bounded subsets of $\kappa$ by $J^{\rm bd}_\kappa$. Given a product of regular cardinals, we denote its true cofinality by ${\rm tcf}$.

We adopt the Jerusalem notation in forcing notions, namely $p\leq q$ means that the condition $q$ gives more information than the condition $p$. We shall use Mathias forcing, relativized to some ultrafilter, and we assume throughout the paper that every ultrafilter is uniform (hence, in particular, non-principal).

\medskip

We thank the referee for many comments, mathematical corrections and a meaningful improvement of the exposition.

\newpage 

\section{Cardinal invariants}
Let us begin with basic definitions of some cardinal invariants. We introduce the general definition, applied to every infinite cardinal $\lambda$ (but in most cases, the definition makes sense only for regular cardinals). Omitting the subscript means that $\lambda=\aleph_0$. Here is the first definition:

\begin{definition}
\label{ssss}
The splitting number $\mathfrak{s}_\lambda$.
\begin{enumerate}
\item [$(\aleph)$] Suppose $B\in[\lambda]^\lambda$ and $S\subseteq\lambda$. $S$ splits $B$ if $|S\cap B|=|(\lambda\setminus S)\cap B|=\lambda$.
\item [$(\beth)$] $\{S_\alpha:\alpha<\kappa\}$ is a splitting family in $\lambda$ if for every $B\in[\lambda]^\lambda$ there exists an ordinal $\alpha<\kappa$ so that $S_\alpha$ splits $B$.
\item [$(\gimel)$] The splitting number $\mathfrak{s}_\lambda$ is the minimal cardinality of a splitting family in $\lambda$.
\end{enumerate}
\end{definition}

The following claim is explicit in \cite{MR2927607} only for the case $\lambda=\aleph_0$ (by our convention, the splitting number is denoted by $\mathfrak{s}$ in this case). Claim 1.3 of \cite{gash962} is also related (but deals with a variant of $\mathfrak{s}$, called the strong splitting number).
For completeness, we repeat the proof here, this time in the general context of $\mathfrak{s}_\lambda$. Notice that the assumption $\lambda<\mathfrak{s}_\lambda$ in the following claim implies that $\lambda$ is weakly compact (we consider $\aleph_0$ as a weakly compact cardinal). A proof appears in \cite{MR1450512} for the case $\lambda$ is regular. We do not know what happens when $\lambda$ is singular (although in some cases a similar result can be proved).

\begin{claim}
\label{ddownward}
The downward positive relation. \newline 
Suppose $\lambda=\cf(\lambda)<\mu<\mathfrak{s}_\lambda$. \newline 
\then\ $\binom{\mu}{\lambda} \rightarrow \binom{\mu}{\lambda}^{1,1}_2$ iff $\cf(\mu)\neq\lambda$.
\end{claim}

\par\noindent\emph{Proof}. \newline 
Assume $\cf(\mu)\neq\lambda$. Let $c:\mu\times\lambda\rightarrow 2$ be any coloring. Set $S_\alpha=\{\gamma\in\lambda:c(\alpha,\gamma)=0\}$ for every $\alpha<\mu$. We collect these sets into the family $\mathcal{F}=\{S_\alpha:\alpha<\mu\}$. Since $|\mathcal{F}|\leq\mu<\mathfrak{s}_\lambda$ we infer that $\mathcal{F}$ is not a splitting family.

Let $B\in[\lambda]^\lambda$ exemplify this fact. It means that $B\subseteq^* S_\alpha$ or $B\subseteq^* (\lambda\setminus S_\alpha)$ for every $\alpha<\mu$. At least one of these options occurs $\mu$-many times, so without loss of generality $B\subseteq^* S_\alpha$ for every $\alpha<\mu$. By the very definition of almost inclusion, for every $\alpha<\mu$ there exists $\beta_\alpha<\lambda$ such that $B\setminus\beta_\alpha\subseteq S_\alpha$ (here we use the regularity of $\lambda$). Since $\cf(\mu)\neq\lambda$ there exists $\beta<\lambda$, and $H_0\in[\mu]^\mu$ so that $\beta_\alpha\leq\beta$ for every $\alpha\in H_0$. 

Let $H_1$ be $B\setminus\beta$, so $H_1\in[\lambda]^\lambda$. Suppose $\alpha\in H_0, \gamma\in H_1$. By the definition of $H_1$, $\gamma\in B\setminus\beta= B\setminus\beta_\alpha$, and since $\alpha\in H_0$ we conclude that $c(\alpha,\gamma)=0$, completing this direction.

Now assume that $\cf(\mu)=\lambda$. Choose a disjoint decomposition $\{A_\gamma:\gamma<\lambda\}$ of $\mu$, such that $|A_\gamma|<\mu$ for every $\gamma<\lambda$. Without loss of generality, the union of every subcollection of less than $\lambda$-many $A_\gamma$-s has size less than $\lambda$. Here we use the assumption $\cf(\mu)=\lambda$.
For every $\alpha<\mu$ let $\xi(\alpha)$ be the unique ordinal so that $\alpha\in A_{\xi(\alpha)}$.
Define $c:\mu\times\lambda\rightarrow 2$ as follows. For $\alpha\in\mu\wedge\beta\in\lambda$ let:

$$
c(\alpha,\beta)=0 \Leftrightarrow \beta\leq\xi(\alpha)
$$

We claim that $c$ exemplifies our claim. Indeed, assume that $|H_0|=\mu$ and $|H_1|=\lambda$. Choose $(\alpha,\beta)\in H_0\times H_1$, and suppose $c(\alpha,\beta)=0$. It means that $\beta\leq\xi(\alpha)$. But $\xi(\alpha)$ is an ordinal below $\lambda$, and $H_1$ is unbounded in $\lambda$, hence one can pick an ordinal $\beta'\in H_1$ so that $\beta'>\xi(\alpha)$. It follows that $c(\alpha,\beta')=1$, so the product $H_0\times H_1$ is not monochromatic in this case. Now suppose $c(\alpha,\beta)=1$. It means that $\xi(\alpha)<\beta$. Clearly, there is some $\alpha'\in H_0$ so that $\xi(\alpha')\geq\beta$. Consequently, $c(\alpha',\beta)=0$, so again $H_0\times H_1$ is not monochromatic, and the proof is completed.

\hfill \qedref{ddownward}

For the next claim we need the following definition:

\begin{definition}
\label{rrrr} The reaping number. \newline 
Let $\lambda$ be an infinite cardinal.
\begin{enumerate}
\item [$(\aleph)$] $\{T_\alpha:\alpha<\kappa\}$ is an unreaped family if there is no $S\in[\lambda]^\lambda$ so that $S$ splits $T_\alpha$ for every $\alpha<\kappa$.
\item [$(\beth)$] the reaping number $\mathfrak{r}_\lambda$ is the minimal cardinality of an unreaped family.
\end{enumerate}
\end{definition}

Our second claim works in the opposite direction to the first claim:

\begin{claim}
\label{uupward} 
The upward positive relation. \newline 
Suppose $\mathfrak{r}_\lambda<\mu\leq 2^\lambda$, $\lambda$ is a regular cardinal. \newline 
\then\ $\binom{\mu}{\lambda} \rightarrow \binom{\mu}{\lambda}^{1,1}_2$ whenever $\cf(\mu)>\mathfrak{r}_\lambda$.
\end{claim}

\par \noindent \emph{Proof}. \newline 
Let $\mathcal{A}\subseteq[\lambda]^\lambda$ exemplify $\mathfrak{r}_\lambda$. It means that $|\mathcal{A}|=\mathfrak{r}_\lambda$, and there is no single $B\in[\lambda]^\lambda$ which splits all the members of $\mathcal{A}$.

Assume $c:\mu\times\lambda\rightarrow 2$ is any coloring. For every $\alpha<\mu$ let $B_\alpha=\{\beta<\lambda:c(\alpha,\beta)=0\}$. Choose $A_\alpha\in\mathcal{A}$ such that $A_\alpha\subseteq^*B_\alpha$ or $A_\alpha\subseteq^*\lambda\setminus B_\alpha$. Without loss of generality, $A_\alpha\subseteq^*B_\alpha$ for every $\alpha<\mu$, so one can choose an ordinal $\beta_\alpha<\lambda$ so that $A_\alpha\setminus\beta_\alpha\subseteq B_\alpha$.

As $\cf(\mu)>\mathfrak{r}_\lambda$, there are $H\in[\mu]^\mu,\beta<\lambda$ and $A\in\mathcal{A}$ such that $\alpha\in H\Rightarrow\beta_\alpha=\beta$ and $A_\alpha= A$. It follows that $c\upharpoonright(H\times A\setminus\beta)=0$, so the proof is completed.

\hfill \qedref{uupward}

Combining the above claims, we can prove the main theorem of this section:

\begin{theorem}
\label{mt}
The main theorem. \newline 
It is consistent that $\binom{\theta}{\omega}\rightarrow\binom{\theta}{\omega}^{1,1}_2$ for every $\aleph_1\leq\theta\leq\ 2^{\aleph_0}$.
\end{theorem}

\par \noindent \emph{Proof}. \newline 
In the model of \cite{MR879489} we have $\mathfrak{r}=\mathfrak{u}=\aleph_1$, while $\mathfrak{s}=\mathfrak{c}=\aleph_2$. By claim \ref{ddownward} we conclude that $\binom{\aleph_1}{\aleph_0}\rightarrow\binom{\aleph_1}{\aleph_0}^{1,1}_2$, and by virtue of claim \ref{uupward} we have $\binom{\aleph_2}{\aleph_0}\rightarrow\binom{\aleph_2}{\aleph_0}^{1,1}_2$, so we are done.

\hfill \qedref{mt}

\begin{corollary}
\label{everyinv}
Polarized relations and cardinal invariants. \newline 
Let $inv$ be any cardinal invariant of the continuum. \newline 
\then\ $\binom{inv}{\omega}\rightarrow\binom{inv}{\omega}^{1,1}_2$ is consistent with ZFC.
\end{corollary}

\hfill \qedref{everyinv}

Notice that $\mathfrak{c}=\aleph_2$ in the model of \cite{MR879489}. Dealing with the dominating number $\mathfrak{d}$, the model in \cite{MR1005010} supplies $\mathfrak{u}=\mu_0<\mu_1=\mathfrak{d}$ for every pair of regular cardinals $(\mu_0,\mu_1)$. It follows that $\binom{\mathfrak{d}}{\omega} \rightarrow\binom{\mathfrak{d}}{\omega}^{1,1}_2$ is consistent for arbitrarily large value of $\mathfrak{d}$, as $\mathfrak{r}\leq\mathfrak{u}$.
We conclude with another open problem from \cite{GaThesis}:

\begin{question}
\label{pppp} The splitting number and the pseudointersection number.
\begin{enumerate}
\item [$(a)$] Is it consistent that $\mathfrak{p}=\mathfrak{s}$ and $\binom{\mathfrak{p}}{\omega}\rightarrow\binom{\mathfrak{p}}{\omega}^{1,1}_2$?
\item [$(b)$] Is it consistent that $\mathfrak{c}=\mathfrak{s}>\aleph_2$ and $\binom{\mathfrak{s}}{\omega}\rightarrow\binom{\mathfrak{s}}{\omega}^{1,1}_2$ (hence $\binom{\theta}{\omega}\rightarrow\binom{\theta}{\omega}^{1,1}_2$ whenever $\cf(\theta)>\aleph_0$)?
\end{enumerate}
\end{question}

Notice that in the above models we have $\mathfrak{p}<\mathfrak{s}$, so a different method is required for this problem. Nevertheless, we believe that a positive answer is consistent for both parts of the question.

\newpage 

\section{Large cardinals}

In this section we deal with uncountable cardinals, with respect to the problems in the previous section. As can be seen, we need some large cardinals assumption. We distinguish two cases. In the first one, $\lambda$ is a regular cardinal. In this case we shall assume that $\lambda$ is a supercompact cardinal, aiming to show that many polarized relations are consistent, above $\lambda$. Secondly, we deal with a singular cardinal.

Let us begin with the regular case. We shall make use of the Mathias forcing, generalized for uncountable cardinals. Notice that for the combinatorial theorems we need a specific version of the Mathias forcing, relativized to some ultrafilter. We begin with the definition of this forcing notion:

\begin{definition}
\label{ggeneralmat}
The generalized Mathias forcing. \newline 
Let $\lambda$ be a supercompact (or even just measurable) cardinal, and $D$ a nonprincipal $\lambda$-complete ultrafilter on $\lambda$. The forcing notion $\mathbb{M}_D^\lambda$ consists of pairs $(a, A)$ such that $a \in [\lambda]^{< \lambda}, A \in D$. For the order, $(a_1, A_1) \leq (a_2, A_2)$ iff $a_1 \subseteq a_2, A_1 \supseteq A_2$ and $a_2 \setminus a_1 \subseteq A_1$.
\end{definition}

Notice that $\mathbb{M}_D^\lambda$ is $\lambda^+$-centered as always $(a,A_1) \parallel (a,A_2)$ and there are only $\lambda=\lambda^{<\lambda}$ many $a$-s. It follows that $\mathbb{M}_D^\lambda$ is $\lambda^+$-cc. Also, $\mathbb{M}_D^\lambda$ is $<\lambda$-directed closed (here we employ the $\lambda$-completeness of the ultrafilter $D$). We emphasize that these properties are preserved by $<\lambda$-support iterations, hence such an iteration collapses no cardinals.

If $\mathbb{M}_D^\lambda$ is a $\lambda$-Mathias forcing, then for defining the Mathias $\lambda$-real we take a generic $G \subseteq \mathbb{M}_D^\lambda$, and define $x_G = \bigcup \{ a : (\exists A \in D)((a, A) \in G) \}$. As in the original Mathias forcing, $x_G$ is endowed with the property $x_G \subseteq^* A \vee x_G \subseteq^* \lambda \setminus A$ for every $A \in [\lambda]^\lambda$ of the ground model. Let us mention another cardinal invariant:

\begin{definition}
\label{uuuu}
The ultrafilter number $\mathfrak{u}_\lambda$. \newline 
Let $\lambda$ be a regular cardinal, and $\mathcal{F}$ a filter on $\lambda$.
\begin{enumerate}
\item [$(\aleph)$] A base $\mathcal{A}$ for $\mathcal{F}$ is a subfamily of $\mathcal{F}$ such that for every $X\in\mathcal{F}$ there is some $Y\in\mathcal{A}$ with the property $Y\subseteq^*X$.
\item [$(\beth)$] $\mathfrak{u}_\lambda$ is the minimal cardinality of a filter base, for some uniform ultrafilter on $\lambda$.
\end{enumerate}
\end{definition}

One can show that $\mathfrak{u}_\lambda>\lambda$ for every $\lambda$.
The following claim employs known facts, so we give just an outline of the proof:

\begin{claim}
\label{uuuandsss} Polarized relations above a supercompact cardinal. \newline 
Suppose $\lambda$ is a supercompact cardinal.
\begin{enumerate}
\item [$(a)$] For every $\mu=\cf(\mu)\in[\lambda^+,2^\lambda]$, one can force $\mathfrak{s}_\lambda=\mu$ without changing the value of $2^\lambda$.
\item [$(b)$] One can force $\mathfrak{u}_\lambda=\lambda^+$ while $2^\lambda$ is arbitrarily large.
\end{enumerate}
\end{claim}

\par\noindent\emph{Outline of proof}.\newline 
For $(a)$ we iterate $\mathbb{M}_D^\lambda$, the length of the iteration being $\mu$. We assume without loss of generality that $\lambda$ is Laver-indestructible, so in particular it remains supercompact (hence measurable) along the iteration. It enables us to choose a $\lambda$-complete ultrafilter at any stage, hence the forcing does not collapse cardinals. We use $<\lambda$-support. It follows that $\mathfrak{s}_\lambda$ equals $\mu$ in the forcing extension. For a detailed proof see also \cite{MR2927607}.

For $(b)$ we use an iteration of length $\lambda^+$. But we choose the $\lambda$-complete ultrafilter (at every stage) more carefully. Along the iteration we create a $\subseteq^*$-decreasing sequence of subsets of $\lambda$. This is done by choosing an ultrafilter which contains the sequence from the previous stages. For the limit stages of the iteration, one has to employ the arguments in \cite{MR1976583}. The main point there is using some prediction principle on $\lambda^+$ in order to make sure that an appropriate ultrafilter is chosen enough times.
At the end, we can show that the sequence (of length $\lambda^+$) generates an ultrafilter, hence $\mathfrak{u}_\lambda=\lambda^+$. We also refer the reader to \cite{adbt} for a detailed proof of this assertion.

\hfill \qedref{uuuandsss}

We indicate that the consistency of $\mathfrak{u}_\lambda=\lambda^+$ while $2^\lambda$ is arbitrarily large is proved for some singular cardinal $\lambda$ in \cite{MR2992547} (but here we deal with a supercompact cardinal).
Let us phrase the following conclusion from the above claim:

\begin{conclusion}
\label{cconc} Many positive relations above a supercompact. \newline 
Suppose $\lambda$ is a supercompact cardinal.
\begin{enumerate}
\item [$(a)$] the positive relation $\binom{\mu}{\lambda}\rightarrow\binom{\mu}{\lambda}^{1,1}_2$ is consistent simultaneously for every regular $\mu$ above $\lambda$ but $2^\lambda$.
\item [$(b)$] the positive relation $\binom{\mu}{\lambda}\rightarrow\binom{\mu}{\lambda}^{1,1}_2$ is consistent simultaneously for every regular $\mu$ in the interval $(\lambda^+,2^\lambda]$.
\end{enumerate}
\end{conclusion}

\par \noindent \emph{Proof}. \newline 
$(a)$ is valid when $\mathfrak{s}_\lambda=2^\lambda$ and $(b)$ holds in a model of $\mathfrak{r}_\lambda=\lambda^+$ (notice that for getting merely $\mathfrak{r}_\lambda=\lambda^+$ we do not need the arguments of \cite{MR1976583}).

\hfill \qedref{cconc}

\begin{question}
\label{uulessthanss} Is it consistent that $\mathfrak{u}_\lambda<\mathfrak{s}_\lambda$ (or at least $\mathfrak{r}_\lambda<\mathfrak{s}_\lambda$) for some uncountable cardinal $\lambda$?
\end{question}

We turn now to the main theorem of this section. We show that getting a positive polarized relation for many cardinals in the interval $(\mu,2^\mu]$ is consistent for some singular cardinal $\mu$ (under some pcf assumptions). In particular, it holds for $\mu^+$.
We shall prove the following:

\begin{theorem}
\label{mmtt} Polarized relations above a singular cardinal. \newline 
Assume $\kappa=\cf(\mu)<\mu<\lambda$, and $\theta<\kappa$. \newline 
If $\circledast$ holds \then\ $\binom{\lambda}{\mu}\rightarrow \binom{\lambda}{\mu}^{1,1}_\theta$ holds, when $\circledast$ means:
\begin{enumerate}
\item [$(a)$] $2^\kappa<\cf(\lambda)$,
\item [$(b)$] $J^{\rm bd}_\kappa\subseteq J$ is an ideal on $\kappa$,
\item [$(c)$] $\langle\lambda_\varepsilon:\varepsilon<\kappa\rangle$ is an increasing sequence of cardinals which tends to $\mu$,
\item [$(d)$] $2^{\lambda_\varepsilon}=\lambda_\varepsilon^+$ for every $\varepsilon<\kappa$,
\item [$(e)$] $\lambda_\varepsilon$ is strongly inaccessible for every $\varepsilon<\kappa$,
\item [$(f)$] $\Upsilon_\ell= {\rm tcf}(\prod\limits_{\varepsilon<\kappa}\lambda_\varepsilon^{+\ell},<_J)$ is well defined for $\ell\in\{0,1\}$,
\item [$(g)$] $\cf(\lambda)\notin\{\Upsilon_0,\Upsilon_1\}$.
\end{enumerate}
\end{theorem}

\par \noindent \emph{Proof}. \newline 
Suppose a coloring $c:\lambda\times\mu\rightarrow\theta$ is given. For every $\alpha<\lambda,\varepsilon<\kappa,\iota<\theta$ we let $A_{\alpha,\varepsilon,\iota}$ be $\{\gamma<\lambda_\varepsilon:c(\alpha,\gamma)=\iota\}$. Fixing $\alpha$ and $\varepsilon$, we have produced a partition $\{A_{\alpha,\varepsilon,\iota}: \iota<\theta\}$ of $\lambda_\varepsilon$ into a small (i.e., just $\theta$-many) number of sets. Enumerate $\mathcal{P}(\lambda_\varepsilon)$ as $\langle B_{\varepsilon,i}:i<\lambda^+_\varepsilon\rangle$. For every $\alpha<\lambda$ and $\iota<\theta$ we define a function $g_{\alpha,\iota}\in\prod\limits_{\varepsilon<\kappa}\lambda_\varepsilon^+$ as follows:

$$
g_{\alpha,\iota}(\varepsilon)={\rm min} \{i<\lambda^+_\varepsilon:A_{\alpha,\varepsilon,\iota}=B_{\varepsilon,i}\}
$$

Here we have used the assumption that $2^{\lambda_\varepsilon}=\lambda^+_\varepsilon$.
For every $\alpha<\lambda$ let $g_\alpha\in\prod\limits_{\varepsilon<\kappa}\lambda_\varepsilon^+$ be defined by $g_\alpha(\varepsilon)= {\rm sup}\{g_{\alpha,\iota}(\varepsilon): \iota<\theta\}$. Note that $g_\alpha(\varepsilon)$ is well defined since each $\lambda_\varepsilon^+$ is regular (but all we need is $\theta<\cf(\lambda_\varepsilon^+)$, to be used in the sequel).

Recall that $\Upsilon_1={\rm tcf}(\prod\limits_{\varepsilon<\kappa}\lambda_\varepsilon^+,<_J)$ and $\cf(\lambda)\neq\Upsilon_1$, hence there exists a function $g\in\prod\limits_{\varepsilon<\kappa}\lambda_\varepsilon^+$ and a set $S_1$ of size $\lambda$ so that $\alpha\in S_1\Rightarrow g_\alpha<_J g$. We may assume, without loss of generality, that $g(\varepsilon)>\lambda_\varepsilon$ for every $\varepsilon<\kappa$. Denote the set $\{\varepsilon<\kappa: g_\alpha(\varepsilon)<g(\varepsilon)\}$ by $u_\alpha$, for every $\alpha<\lambda$. Since $2^\kappa<\cf(\lambda)$, there are $u\subseteq\kappa$ and $S_2\in[S_1]^\lambda$ such that $u=\kappa\ {\rm mod\ }J$ and $\alpha\in S_2\Rightarrow u_\alpha=u$. Without loss of generality, $u=\kappa$.

Take a closer look at the collection $\{B_{\varepsilon,i}:i<g(\varepsilon)\}$ (for every $\varepsilon<\kappa$). By the nature of the function $g$, this is a family of $\lambda_\varepsilon$-many sets, hence we can enumerate its members as $\{B_{\varepsilon,i}^1:i<\lambda_\varepsilon\}$. Notice that for every $\alpha\in S_2,\varepsilon<\kappa$ and $\iota<\theta$ we know that $A_{\alpha,\varepsilon,\iota}\in \{B_{\varepsilon,i}^1:i<\lambda_\varepsilon\}$.

We need another round of unifying. By the same token as above, we define for every $\alpha\in S_2$ and $\iota<\theta$ the function $h_{\alpha,\iota}\in \prod\limits_{\varepsilon<\kappa}\lambda_\varepsilon$ as follows:

$$
h_{\alpha,\iota}(\varepsilon)={\rm min}\{i<\lambda_\varepsilon: B_{\varepsilon,i}^1=A_{\alpha,\varepsilon,\iota}\}
$$

Now for $\alpha\in S_2$ set $h_\alpha(\varepsilon)={\rm sup}\{h_{\alpha,\iota}(\varepsilon)+1:\iota<\theta\}$ (for every $\varepsilon<\kappa$). Again, by our assumptions, $h_\alpha$ belongs to the product $\prod\limits_{\varepsilon<\kappa}\lambda_\varepsilon$ for every $\alpha\in S_2$.
Since $\cf(\lambda)\neq\Upsilon_0$ (recall that $\Upsilon_0= {\rm tcf}(\prod\limits_{\varepsilon<\kappa}\lambda_\varepsilon,<_J)$) we can choose a function $h$ which bounds many $h_\alpha$-s. In other words, there are $h$ and $S_3\in[S_2]^\lambda$ so that $\alpha\in S_3\Rightarrow h_\alpha<_Jh$.

Let $v_\alpha$ be the set $\{\varepsilon<\kappa: h_\alpha(\varepsilon)<h(\varepsilon)\}$, for every $\alpha\in S_3$. As before, since $2^\kappa<\cf(\lambda)$ one can find $v\subseteq\kappa$ and $S_4\in[S_3]^\lambda$ so that $\alpha\in S_4\Rightarrow v_\alpha=v$. Without loss of generality we assume, as usual, that $v=\kappa$.

For every $\varepsilon<\kappa$ we define an equivalence relation $E_\varepsilon$ on $\lambda_\varepsilon$ as follows:

$$
\forall\gamma_1,\gamma_2\in\lambda_\varepsilon, \gamma_1E_\varepsilon\gamma_2 \Leftrightarrow (\gamma_1\in B_{\varepsilon,j}^1\equiv\gamma_2\in B_{\varepsilon,j}^1,\forall j<h(\varepsilon))
$$

Observe that $E_\varepsilon$ has less than $\lambda_\varepsilon$ equivalence classes for every $\varepsilon<\kappa$, since each $\lambda_\varepsilon$ is an inaccessible cardinal. Consequently, we can choose an equivalence class $X_\varepsilon$ of size $\lambda_\varepsilon$ in each $E_\varepsilon$. For every $\alpha\in S_4$ let $\iota_{\alpha,\varepsilon}<\theta$ be the color associated with $X_\varepsilon$ (i.e., $c(\alpha,\gamma)=\iota_{\alpha,\varepsilon}$ for every $\gamma\in X_\varepsilon$).

We arrived at the last stage of unifying $\varepsilon$-s. For every $\alpha\in S_4$ there is a color $\iota_\alpha$ so that the set $w_\alpha=\{\varepsilon<\kappa: \iota_{\alpha,\varepsilon}=\iota_\alpha\}$ is of size $\kappa$. Hence there are a color $\iota<\theta,w\in[\kappa]^\kappa$ and $S_5\in[S_4]^\lambda$ such that $\alpha\in S_5\Rightarrow \iota_\alpha= \iota,w_\alpha=w$.

Set $A=S_5$ and $B=\bigcup\{X_\varepsilon:\varepsilon\in w\}$. Clearly, $A\in [\lambda]^\lambda,B\in[\mu]^\mu$. We claim that the product $A\times B$ exemplifies the positive relation $\binom{\lambda}{\mu}\rightarrow \binom{\lambda}{\mu}^{1,1}_2$. Indeed, if $\alpha\in A$ and $\beta\in B$ then $\alpha\in S_5$ and $\beta\in X_\varepsilon$ for some $\varepsilon\in w$. Consequently, $\iota_\alpha=\iota$ (for this specific $\alpha$) and $c(\alpha,\beta)=\iota_{\alpha,\varepsilon}=\iota_\alpha=\iota$ so we are done.

\hfill \qedref{mmtt}

\begin{corollary}
\label{muplus} Positive relation for successor of singular. \newline 
Suppose $(\kappa,\mu,\mu^+)$ satisfy $\circledast$ of Theorem \ref{mmtt} (stipulating $\mu^+=\lambda$).
\then\ $\binom{\mu^+}{\mu}\rightarrow \binom{\mu^+}{\mu}^{1,1}_2$. \newline 
In particular, this positive relation is consistent with ZFC.
\end{corollary}

\par \noindent \emph{Proof}. \newline 
We refer to \cite{GaSh949}, where the assumptions of the theorem are forced (and in fact, much more), but see also the discussion following the next remark below.

\hfill \qedref{muplus}

\begin{remark}
\label{949}
A similar result is forced in \cite{GaSh949}, under the assumption that $\mu$ is a singular cardinal which is a limit of measurables. In the forcing extension of \cite{GaSh949}, one has to admit the existence of a supercompact cardinal in the ground model. Nevertheless, the polarized relation there is slightly stronger. Being a limit of measurables entails $\binom{\lambda}{\mu}\rightarrow \binom{\lambda}{\mu}^{1,<\omega}_2$ there (which means that for every $c:\lambda\times[\mu]^{<\omega}\rightarrow 2$ there are $H_0\in[\lambda]^\lambda,H_1\in[\mu]^\mu$ such that for every $n\in\omega, c\upharpoonright(H_0\times[H_1]^n)$ is constant).

We also indicate that the assumption $2^{\lambda_\varepsilon}=\lambda^+_\varepsilon$ is stronger than needed here. The value of $2^{\lambda_\varepsilon}$ can be replaced by a larger cardinal, provided that all the relevant products have true cofinality. Anyhow, some restriction should be imposed. If $2^{\lambda_\varepsilon}= \lambda_\varepsilon^{+\zeta(\varepsilon)}$ for every $\varepsilon<\kappa$, and the sequence $\langle\zeta(\varepsilon):\varepsilon<\kappa\rangle$ tends to $\mu$, then the argument breaks down.
\end{remark}

\hfill \qedref{949}

We can modify the proof above, to include another case. The consistency proof of the assumptions below is similar to those of Theorem \ref{mmtt}, yet comparing to \cite{GaSh949} we need less than supercompactness (in both theorems). A sufficient assumption in order to force the assumptions of these theorems is the existence of a strong cardinal in the ground model (and even slightly less, namely a $\tau$-strong cardinal for some suitable $\tau$). We hope to shed light on this subject in \cite{GaMgShF1211}.

\begin{theorem}
\label{sssinggg} Positive relation for limit of strong limit cardinals. \newline 
Assume $\kappa=\cf(\mu)<\mu<\lambda$, and $\theta<\kappa$. \newline 
If $\circledcirc$ holds \then\ $\binom{\lambda}{\mu}\rightarrow \binom{\lambda}{\mu}^{1,1}_\theta$ holds, when $\circledcirc$ means:
\begin{enumerate}
\item [$(a)$] $2^\kappa<\cf(\lambda)$,
\item [$(b)$] $J^{\rm bd}_\kappa\subseteq J$ is an ideal on $\kappa$,
\item [$(c)$] $\langle\lambda_\varepsilon:\varepsilon<\kappa\rangle$ is an increasing sequence of cardinals which tends to $\mu$,
\item [$(d)$] $2^{\lambda_\varepsilon}=\lambda_\varepsilon^+$ for every $\varepsilon<\kappa$,
\item [$(e)$] $\lambda_\varepsilon$ is strong limit and $\cf(\lambda_\varepsilon)>\kappa$ for every $\varepsilon<\kappa$,
\item [$(f)$] $\prod\limits_{\varepsilon<\kappa}\cf(\lambda_\varepsilon)< \cf(\lambda)$,
\item [$(g)$] $\Upsilon_\ell= {\rm tcf}(\prod\limits_{\varepsilon<\kappa}\lambda_\varepsilon^{+\ell},<_J)$ is well defined for $\ell\in\{0,1\}$,
\item [$(h)$] $\cf(\lambda)\notin\{\Upsilon_0,\Upsilon_1\}$.
\end{enumerate}
\end{theorem}

\par\noindent\emph{Proof}.\newline 
Proceed as in the proof of Theorem \ref{mmtt}, till the stage of defining the equivalence relations $E_\varepsilon$ on each $\lambda_\varepsilon$. At this stage we have isolated a large equivalence class (for every $\varepsilon<\kappa$), using the regularity of $\lambda_\varepsilon$. But here, $\lambda_\varepsilon$ is a singular cardinal, so we have to be more careful.

For every $\varepsilon<\kappa$ we choose a sequence $\langle X_{\varepsilon,j} :j<\cf(\lambda_\varepsilon)\rangle$ so that each $X_{\varepsilon,j}$ is an equivalence class of $E_\varepsilon$, and $\Sigma\{|X_{\varepsilon,j}|: j<\cf(\lambda_\varepsilon)\}=\lambda_\varepsilon$.
For every $\alpha\in S_4, \varepsilon<\kappa$ and $j<\cf(\lambda_\varepsilon)$ we choose a color $\iota_{\alpha,\varepsilon,j}<\theta$ so that:

$$
\gamma\in X_{\varepsilon,j}\Rightarrow c(\alpha,\gamma)=\iota_{\alpha,\varepsilon,j}
$$

We claim that there are $S_5\in[S_4]^\lambda$ and a sequence of colors $\langle\iota_{\varepsilon,j}:\varepsilon<\kappa,j< \cf(\lambda_\varepsilon)\rangle$ such that $\alpha\in S_5\Rightarrow \iota_{\alpha,\varepsilon,j}=\iota_{\varepsilon,j}$ (here we use assumption $(f)$ of the present theorem). Moreover, there is a single color $\iota<\theta$ so that $\Sigma\{|X_{\varepsilon,j}|: \iota_{\varepsilon,j}=\iota, \varepsilon<\kappa,j<\cf(\lambda_\varepsilon)\}=\mu$.
For this, notice that $\mu=\Sigma_{\varepsilon<\kappa} \lambda_\varepsilon = \Sigma_{\varepsilon<\kappa}\Sigma \{|X_{\varepsilon,j}|: j<\cf(\lambda_\varepsilon)\}$.

Now we can define $A=S_5$ and $B=\bigcup \{X_{\varepsilon,j}: \iota_{\varepsilon,j}=\iota, \varepsilon<\kappa, j<\cf(\lambda_\varepsilon)\}$. It follows that $A\in[\lambda]^\lambda$ and $B\in[\mu]^\mu$. Since the product $A\times B$ is monochromatic, we are done.

\hfill \qedref{sssinggg}

\begin{remark}
\label{ffff} The assumption $\prod\limits_{\varepsilon<\kappa}\cf(\lambda_\varepsilon)< \cf(\lambda)$ ((f) in the last theorem) can be omitted. We have to choose an equivalence class $X_\varepsilon$ of $E_\varepsilon$ of size at least $(\sum\limits_{\zeta<\varepsilon}\lambda_\zeta)^+$ so that ${\rm min}(X_\varepsilon)>\sum\{\lambda_\zeta:\zeta<\varepsilon\}$. But in some sense we get less.
\end{remark}

We conclude this section with the following:

\begin{proposition}
\label{alll} It is consistent that there is a singular $\mu,\kappa=\cf(\mu)$, such that $\binom{\lambda}{\mu}\rightarrow\binom{\lambda}{\mu}^{1,1}_2$ holds for all $\lambda\in(\mu,2^\mu]$.
\end{proposition}

\par\noindent\emph{Proof}.\newline 
By enlarging $2^\mu$ to a large enough value below $\mu^{+\omega}$, one can choose two sequences, $\langle\lambda_\varepsilon :\varepsilon<\kappa\rangle$ and $\langle\kappa_\varepsilon :\varepsilon<\kappa\rangle$ of inaccessibles (for simplicity) whose limit is $\mu$, and $\{\Upsilon_0^{\bar{\lambda}},\Upsilon_1^{\bar{\lambda}}\} \cap \{\Upsilon_0^{\bar{\kappa}},\Upsilon_1^{\bar{\kappa}}\}=\emptyset$ (the ideal we use is $J_\kappa^{\rm bd}$, see, for instance, the models constructed in \cite{MR1900900}). 
Notice that all the cardinals in the interval $[\mu^+,2^\mu]$ are regular.
Now use Theorem \ref{mmtt}.

\hfill \qedref{alll}

\newpage 

\bibliographystyle{amsplain}
\bibliography{arlist}

\end{document}